\magnification=\magstep1
\hsize16truecm
\vsize23.5truecm
\topskip=1truecm
\raggedbottom
\abovedisplayskip=3mm
\belowdisplayskip=3mm
\abovedisplayshortskip=0mm
\belowdisplayshortskip=2mm
\normalbaselineskip=12pt
\normalbaselines
\font\titlefont= cmcsc10 at 12pt
\def\F{\Bbb F}
\def\R{\Bbb R}

\def\Z{\Bbb Z}
\def\Q{\Bbb Q}
\def\P{\Bbb P}

\def\G{\Bbb G}

\def\E{\Bbb E}

\def\H{\Bbb H}

\def\lllongrightarrow{\relbar\joinrel\relbar\joinrel\relbar\joinrel\rightarrow}
%
%
\catcode`\@=11
\font\tenmsa=msam10
\font\sevenmsa=msam7
\font\fivemsa=msam5
\font\tenmsb=msbm10
\font\sevenmsb=msbm7
\font\fivemsb=msbm5
\newfam\msafam
\newfam\msbfam
\textfont\msafam=\tenmsa  \scriptfont\msafam=\sevenmsa
  \scriptscriptfont\msafam=\fivemsa
\textfont\msbfam=\tenmsb  \scriptfont\msbfam=\sevenmsb
  \scriptscriptfont\msbfam=\fivemsb
\def\hexnumber@#1{\ifcase#1 0\or1\or2\or3\or4\or5\or6\or7\or8\or9\or
	A\or B\or C\or D\or E\or F\fi }
\edef\msa@{\hexnumber@\msafam}
\edef\msb@{\hexnumber@\msbfam}
\mathchardef\square="0\msa@03
\mathchardef\subsetneq="3\msb@28
\mathchardef\ltimes="2\msb@6E
\mathchardef\rtimes="2\msb@6F
\def\Bbb{\ifmmode\let\next\Bbb@\else
	\def\next{\errmessage{Use \string\Bbb\space only in math mode}}\fi\next}
\def\Bbb@#1{{\Bbb@@{#1}}}
\def\Bbb@@#1{\fam\msbfam#1}
\catcode`\@=12
%
%
%
\vskip 6.5pc
\noindent
\font\eighteenbf=cmbx10 scaled\magstep3
\vskip 2.0pc
\centerline{\eighteenbf  Cycles on the Moduli Space }
\bigskip
\centerline{\eighteenbf of Abelian Varieties}\footnote{}{alg-geom 9605011, 
revised version: June 17, 1996}
\noindent
\vskip 2pc
\font\titlefont=cmcsc10 at 12pt
\centerline{\titlefont Gerard van der Geer}
\vskip 2.0pc
\noindent
\bigskip
\noindent
\centerline{\bf  Introduction}
\smallskip
In this paper I present a number of results on cycles on the moduli
space ${\cal A}_g$ of principally polarized abelian varieties of
dimension $g$. The results  on the tautological ring are my own
work, the results on the torsion of
$\lambda_g$ and on the cycle classes of the Ekedahl-Oort stratification
are joint work with Torsten Ekedahl and some of the results on curves are
joint work with Carel Faber. Our results include:
\smallskip
\noindent
\item{$\bullet$} A description of the tautological subring -- generated by  
the Chern classes $\lambda_i$ of the Hodge bundle $\E$ -- of the Chow
ring of ${\cal A}_g$ and its compactifications.\par 
\smallskip
\noindent
\item{$\bullet$} As a corollary we find the Hirzebruch-Mumford proportionality
theorem for ${\cal A}_g$ and bounds for the dimension of complete
subvarieties of ${\cal A}_g$. \par
\smallskip
\noindent
\item{$\bullet$} A bound for the order of the torsion of the top Chern class
$\lambda_g$ of the Hodge bundle.\par
\smallskip
\noindent
\item{$\bullet$} A description of the Ekedahl-Oort stratification of ${\cal
A}_g \otimes \F_p$ in terms of degeneracy loci of a map between flag
bundles.\par
\smallskip
\noindent
\item{$\bullet$}  The description of the Chow classes of the strata of this
stratification. This includes as special cases formulas for the
classes of loci like $p$-rank $\leq f$ locus or $a$-number $\geq a$
locus. \par 
\smallskip
\noindent
\item{$\bullet$}  The irreducibility of the locus $T_a$  of abelian varieties of
$a$-number $\geq a$ for $a<g$. \par
\smallskip
\noindent
\item{$\bullet$}  A computation of this stratification for hyperelliptic curves
of $2$-rank $0$ in characteristic~2. \par
\smallskip
\noindent
\item{$\bullet$} A formula for the class of the supersingular locus for low
genera.
\smallskip
It is a pleasure for me to acknowledge pleasant cooperation with and help
from Torsten Ekedahl and Carel Faber and useful discussions with
H\'el\`ene Esnault, Frans Oort and Piotr Pragacz. I also would like to
thank Kenji Ueno for inviting me to Kyoto where I found the time to write
this paper.
\bigskip
\smallskip\noindent
{\bf \S 1. The Tautological Subring of} ${\cal A}_g$ {\sl  and of}
${\tilde {\cal A}}_g$. 
\bigskip
\noindent
Let ${\cal A}_g/ \Z$ denote the moduli stack of principally polarized
abelian varieties of dimension~$g$. This is an irreducible algebraic
stack of relative dimension $g(g+1)/2$. All the fibres ${\cal
A}_g\otimes \F_p$ and ${\cal A}_g \otimes \Q$ are irreducible. This stack 
carries a locally free sheaf $\E$ of rank $g$ (``the Hodge bundle")
defined by giving for every morphism $S \to {\cal A}_g$ a locally free sheaf of
rank $g$ which is compatible with pull backs. It is defined as
$s^*\Omega_{A/S}$, where $A/S$ is the principally polarized abelian variety
corresponding to  $S \to {\cal A}_g$ and $s$
is the zero section of $A/S$. If $\pi: {\cal X}_g \to {\cal A}_g$ is the
universal abelian variety we have $\Omega_{ {\cal X}_g/ {\cal A}_g}=
\pi^*(\E)$.
\smallskip
Let $\tilde {\cal A}_g$ be a smooth toroidal compactification of ${\cal A}_g$.
The Hodge bundle can be extended to a locally free sheaf on $\tilde {\cal A}_g$.
\smallskip
The Chern classes $\lambda_i$ of the Hodge bundle $\E$ are
defined over
$\Z$ and give rise to classes $\lambda_i$ in $CH^*({\cal A}_g)$, and in 
$CH^*(\tilde {\cal A}_g)$. They generate  subrings ($\Q$-subalgebras) of
$CH_{\Q}^*({\cal A}_g)$  and of $CH_{\Q}^*(\tilde {\cal A}_g)$ which we 
shall call the {\sl tautological subrings}.  
\smallskip
We shall first describe these tautological  subrings. It  will
turn out that  the tautological subring of $CH_{\Q}^*({\cal A}_g)$ is
isomorphic to the cohomology ring $R_{g-1}$  of the compact dual of the
Siegel upper half space of degree
$g-1$, while the tautological subring of $CH_{\Q}^*(\tilde {\cal A}_g)$ is
isomorphic to $R_g$. This cohomology ring  is of the form 
$$
R_g= \Q[u_1,\ldots, u_g]/((1+u_1+\ldots +u_g)(1-u_1+u_2-\ldots + (-1)^gu_g)-1)
$$
and is a Gorenstein ring. As a corollary of this we find  the Proportionality
Principle of Hirzebruch and Mumford for ${\cal A}_g$ and its compactifications. 
\bigskip
We have the following  relation for the Chern classes
$\lambda_i$ on ${\cal A}_g$.
\smallskip
\noindent
\proclaim (1.1) Theorem. The Chern classes $\lambda_i$ in $CH_{\Q}^*( {\cal
A}_g)$ satisfy the relation
$$
(1+\lambda_1+\ldots +\lambda_g)(1-\lambda_1+\ldots +(-1)^g
\lambda_g)=1. \eqno(1)
$$
\par
The idea of the proof is to apply the Grothendieck-Riemann-Roch theorem to
the theta divisor on the universal abelian variety ${\cal X}_g$ over
${\cal A}_g$. We choose this divisor (on a level cover) so
that its restriction $s^*(\Theta)$, with $s$ the zero section, is trivial on
${\cal A}_g$ and apply Grothendieck-Riemann-Roch to the line bundle
$L=O(\Theta)$: 
$$ \eqalign
	{ 
ch(\pi_!L)&= \pi_*(ch(L) \cdot  {\rm Td}(\Omega_{{\cal X}_g/ {\cal
A}_g}^1)^{\vee})\cr
	&=\pi_*( ch(L) \cdot {\rm Td}(\pi^*(\E^{\vee})))\cr
	&=\pi_*(ch(L)) \cdot  {\rm Td}(\E^{\vee})\cr
	}
$$
by the projection formula. Since $R^i\pi_*(L) =0$ for $i>0$ it follows
that $\pi_!(L)$ is a vector bundle and it is of rank one since
$\Theta$ is a principal polarization. We write $c_1(\pi_!(L))= \theta$
and find
$$
\sum_{k=0}^{\infty} {\theta^k\over k!} = \pi_*(\sum_{k=0}^{\infty}
{\Theta^{g+k}\over (g+k)!}) \cdot {\rm Td}(\E^{\vee}). \eqno(2)
$$
Comparison of the term of degree $1$ gives 
$$
\theta = - \lambda_1/2 + \pi_*(\Theta^{g+1})/(g+1)! .
$$
Replace now $L$ by $L^{\otimes n}$. The term of degree $k$  in $\pi_*(\sum
\Theta^{g+k}/ (g+k)!)$ changes by a factor $n^{g+k}$. But $\pi_!(L^{n})$
is a numerical function of degree $\leq n^g$,  cf.\ the
arguments of Chai-Faltings on p.\ 26 of [F-C]. Or, alternatively, using
the Heisenberg group we find  on a suitable cover of ${\cal A}_g$ 
$$
\pi_!L^n = \pi_!L\otimes A
$$
where $A$ is the standard irreducible representation of dimension $n^g$ of the
Heisenberg group. Therefore, up to torsion we have
$$
ch(\pi_{!} L^n) = n^g ch(\pi_{!}L).
$$
This then proves by induction that
$$ \pi_*(\sum_{k=0}^{\infty} {\Theta^{g+k}\over (g+k)!}) = 1 \in
CH_{\Q}^*({\cal A}_g). \eqno(3) 
$$   
In particular we get 
$$
2\theta = -\lambda_1\qquad {\hbox{\rm (the ``key formula")}}.
$$
Therefore, if we write $\lambda_j$ as the $j$th symmetric function of
$\alpha_1,\ldots, \alpha_g$ and if we use ${\rm Td}(\E^{\vee}) = \prod
(\alpha_i /(e^{\alpha_i} -1)$) the identity (2) becomes  $$
\prod_{i=1}^g  { e^{\alpha_i /2} - e^{-\alpha_i/2} \over \alpha_i} =1.
$$
and implies ${\rm Td}(\E\oplus\E^{\vee})=1$. This is equivalent with
$ch_{2k}(\E)=0$ for $k\geq 1$  and with~(1).
$\square$
\smallskip
\noindent
\smallskip
\noindent
\proclaim (1.2) Proposition. In $CH_{\Q}^*({\cal A}_{g})$ we have $\lambda_g=0$.
\par
\noindent 
{\sl Proof.} We apply GRR to the structure sheaf $O_{\cal X}$ of the universal
abelian variety ${\cal X}$ over the stack ${\cal A}_g$. We get
$$
\eqalign{ ch(\pi_{!} O_{\cal X}) &= \pi_*(ch(O_{\cal X})\cdot {\rm
Td}(\Omega^1)^{\vee})
\cr
&=\pi_*(1) {\rm Td}(\E^{\vee}), \cr }
$$
which gives
$$
ch(1-\E^{\vee} + \wedge^2 \E^{\vee} - \ldots +(-1)^g \wedge^g \E^{\vee}) 
= \pi_* (1) {\rm Td}(\E)^{\vee}=0.
$$
Now for a vector bundle $B$ of rank $r$ we have in general the relation (see
[B-S])
$$
\sum_{j=0}^r (-1)^j ch(\wedge^j B^{\vee}) = c_r(B) {\rm Td}(B)^{-1}.
$$
This gives: $\lambda_g=0$ in $CH^*({\cal A}_g)$. $\square$
\bigskip
\noindent
The ring $R_g$ is the quotient of the graded ring $\Q[u_1,\ldots,u_g]$ with
generators $u_i$ of degree $i$  by the relation
$$
(1+u_1+u_2+\ldots+u_g)(1-u_1+u_2- \ldots +(-1)^gu_g)-1=0.
$$
This relation implies (by induction)
$$
u_gu_{g-1}\cdots u_{k+1}u_k^2=0 \qquad {\rm for }\quad k=1,\ldots,g.\eqno(4)
$$
This ring has additive generators
$$
u_1^{\epsilon_1}u_2^{\epsilon_2}\ldots u_g^{\epsilon_g} \qquad {\rm with}
\quad \epsilon_j \in \{ 0,1\}.
$$
Obviously, we have $R_g/(u_g) \cong R_{g-1}$.
\smallskip
It follows that the tautological subring of
$CH_{\Q}^*({\cal A}_g)$ is a homomorphic image of the ring $R_g/(u_g) \cong
R_{g-1}$. 
\bigskip
Now consider the moduli space ${\cal A}_g \otimes \F_p$. It contains the loci
$V_f= V_f(p)$ of abelian varieties with $p$-rank $\leq f$. Their closures in
${\cal A}_g^*$ and in $\tilde {\cal A}_g\otimes \F_p$ define loci again denoted
$V_f$. These loci have been studied by Oort and Norman (cf. [O1 , N-O]).
\smallskip
\noindent
\proclaim (1.3) Lemma. The subvariety $V_f$ is complete in the moduli space
$\tilde{\cal A}_g^{(f)} \subset \tilde {\cal A}_g\otimes \F_p$ of rank$\leq f$
degenerations; in particular $V_0$ is complete in ${\cal A}_g$.
\par 
\smallskip
\noindent
\proclaim (1.4) Corollary. We have $\lambda_1^{g(g-1)/2 +f} \neq 0$ on
$\tilde{\cal A}_g^{(f)}$. \par
\smallskip
\noindent
{\sl Proof.} Observe that $\det (\E)$ is an ample line bundle
(its sections are modular forms) and so $\lambda_1$ is ample on ${\cal
A}_g^*$, see [M-B]. On a complete variety of dimension $d$ the $d$-th
power of an ample divisor is non-zero. $\square$
\smallskip \noindent
\proclaim  (1.5) Theorem. The tautological subring of $CH_{\Q}^*({\cal
A}_{g})$ generated by  the $\lambda_i$ is isomorphic to $R_{g-1}$. \par
\noindent
{\sl Proof.} By the relation $(1+\lambda_1 + \ldots + \lambda_g)
(1-\lambda_1+ \ldots +(-1)^g \lambda_g)=1$ we get a quotient ring of $R_g$. 
Since moreover  $\lambda_g=0$ we get a quotient of $R_g/(u_g)\cong R_{g-1}$.
This ring $R_{g-1}$ is Gorenstein and its top degree
elements are proportional. We now use the fact that ${\cal A}_g\otimes
\F_p$ has a complete subvariety of codimension
$g$, namely the $p$-rank zero locus. The existence of of a complete subvariety of
dimension  $g(g-1)/2$ in ${\cal A}_{g}\otimes \F_p$  for every prime $p$
and the ampleness of $\lambda_1$ on ${\cal A}_g \otimes \F_p$ imply 
$\lambda_1^{g(g-1)/2}\neq 0$. This implies that in the quotient of
$R_{g-1}$ the one-dimensional socle does not map to zero, hence that the
quotient is isomorphic to $R_{g-1}$.  Or more explicitly, consider the set
of $2^{g-1}$ generators of the form
$$
\lambda_1^{\epsilon_1}\lambda_2^{\epsilon_2} \cdots 
\lambda_{g-1}^{\epsilon_{g-1}} 
\quad{\rm with}\quad \epsilon_i\in \{ 0,1\}.
$$
Order these elements $\lambda_{\epsilon}$ with $\epsilon  \in \{ 0,1\}^{g-1}$ 
lexicographically. Suppose we have a relation
$$
\sum a_{\epsilon} \lambda_{\epsilon} = 0 
$$
in $CH_{\Q}^*({\cal A}_g)$. Suppose that $\epsilon'$ is the `smallest'
exponent. Let $\epsilon''$ be the complementary exponent with 
$\epsilon' + \epsilon''=(1,\ldots,1)$. Then we have
$$
\lambda_{\epsilon''} (\sum a_{\epsilon} \lambda_{\epsilon})= a_{\epsilon'}
\lambda_{g-1} \cdots \lambda_1 =0,
$$
and this implies that $a_{\epsilon'} =0$. By induction the theorem follows. 
$\square$
\smallskip
\noindent
\proclaim (1.6) Corollary. In $CH_{\Q}^*({\cal A}_g)$ we have: 
$\lambda_1^{g(g-1)/2} \neq 0$ 
and
$\lambda_1^{1+ g(g-1)/2} = 0$. \par  
\noindent
{\sl Proof.} The first statement was given in (1.4). In the ring $R_{g-1}$ all
top monomials (of degree
$g(g-1)/2$) are pro\-portional, hence $\lambda_1^{g(g-1)/2}$ is a non-zero
multiple of the top monomial  
$\lambda_{g-1}\lambda_{g-2}\cdots \lambda_1$, so
$\lambda_1^{g(g-1)/2 + 1}$ is a non-zero multiple of 
$\lambda_{g-1}\lambda_{g-2}\cdots \lambda_2\lambda_1^2$, which is zero 
by (4).
$\square$
\smallskip
\noindent
\proclaim (1.7) Corollary. Let $\F$ be a field. A complete subvariety of 
 ${\cal A}_g\otimes \F$ has codimension~$\geq g$. \par
\smallskip
\noindent
{\sl Proof.} If $Z$ is a complete subvariety of dimension $m$ then
$\lambda_1^m\neq 0$ on $Z$. Since $g(g-1)/2$ is the highest power of
$\lambda_1$ which is not zero on ${\cal A}_g$ the result follows.
$\square$
\smallskip
For a discussion of questions concerning complete subvarieties of
${\cal A}_g$ we refer to [O2]. 
\bigskip
Now we come to the tautological ring of a suitable toroidal
compactification of ${\cal A}_g$. We shall consider various 
compactifications of ${\cal A}_g$. We choose a suitable compactification
$\tilde{\cal A}_g$ as constructed in [F-C]. We let ${\cal A}_g^*$ be the
minimal (`Satake') compactification as in [F-C]. Furthermore, we let
$\tilde{\cal A}_g^{(1)}$ be the moduli space of rank $1$-degenerations,
i.e. the inverse image of  ${\cal A}_g \cup {\cal A}_{g-1} \subset {\cal
A}_g^*$ under the natural map $ q: \tilde{\cal A}_g \to {\cal A}_g^*$.
This space ${\cal A}_g' = \tilde{\cal A}_g^{(1)}$ does {\sl not} depend on
a choice
$\tilde {\cal A}_g$ of compactification of ${\cal A}_g$. See [M3].  
\smallskip
Let
${\cal G}$ be the `universal' semi-abelian variety over $\tilde{\cal A}_g$
with zero section $s: \tilde{\cal A}_g \to {\cal G}$. Then we have 
$\E= s^*({\rm Lie}{ \cal G})^{\vee}$. 
Moreover, let $D= \tilde{\cal A}_g - {\cal A}_g$ be the divisor at
infinity. Then we have the isomorphism (cf. [F-C]. p.\ 117):
$$
{\rm Sym}^2(\E)\cong \Omega^1(\log D).
$$
Since the class $\lambda_g$ vanishes on ${\cal A}_g$ it can be 
represented in the form $\lambda_g= i_*(x)$ with $x$ a class in
$CH_{\Q}^{g-1}(D)$ and $i\colon D \to \tilde{\cal A}_g$ the embedding
of the boundary. By applying Grothendieck-Riemann-Roch to the structure sheaf
on the semi-abelian variety over $\tilde{\cal A}_g$ (as in (1.2)) 
one gets an  expression for $x$ in $CH_{\Q}^*(D)$, see [EFG].

As we shall see later in characteristic $p$ a multiple of the class
$\lambda_{g}$ is represented in the Chow ring by the class of a complete
subvariety $V_0$ of ${\cal A}_g$. More generally, a multiple of
$\lambda_{g-f}$ is represented by a complete subvariety of ${\cal
A}_{g}^{(f)}$, the moduli space of rank $\leq f$ degenerations.
\smallskip
\noindent
\proclaim (1.8) Lemma.  The class $q_*(\lambda_g)$ in $CH_{\Q}^g({\cal
A}_g^*)$ is  represented by a multiple of the fundamental class of the
boundary $B_g^*={\cal A}_g^* - {\cal A}_g$. \par
\smallskip
Indeed, $\lambda_g$ is zero on ${\cal A}_g$; for dimension
reasons $q_*\lambda_g$ is represented by a multiple of the fundamental
class of $B_g^*$.
\smallskip
\noindent
\proclaim (1.9) Proposition. The cycle class $[B_g^*]$ of
the boundary is the same in the Chow group $CH_{\Q}^g({\cal A}_g^*)$ as a
multiple of the class of the image of ${\cal A}_{g-1}$ in ${\cal A}_g$
under the map
$[X]
\mapsto [X \times E]$, with $E$ a generic elliptic curve.
\par
\smallskip
\noindent
{\sl Proof.} Consider (for $g>2$) the space ${\cal A}_{g-1,1}\sim {\cal
A}_{g-1} \times {\cal A}_1$ in ${\cal A}_g$. Since ${\cal A}_1$ is the
affine $j$-line we find a rational equivalence between the cycle class
of a generic fibre ${\cal A}_{g-1}
\times \{ j \}$ and a multiple of the fundamental class of the boundary
$B_g^*$.  
$\square$
\bigskip
Let $B_g$ be the cycle on ${\cal A}_g^{(1)}$ defined by the semi-abelian
varieties which are trivial extensions
$$
1 \to \G_m \to G \to X_{g-1} \to 0,
$$
where $X_{g-1}$ is a $(g-1)$-dimensional abelian variety. The cycle
$B_g\sim {\cal A}_{g-1}$ is of codimension $g$ in ${\cal A}_g^{(1)}$ and
extends to
$\tilde{\cal A}_g$.
\smallskip
The following proposition describes the class of $B_g$ in {\sl
cohomology}. 
\smallskip
\noindent
\proclaim (1.10) Proposition. The cohomology  class of $B_g$ and of
${\cal A}_{g-1} \times \{j\}$ are both equal to 
$(-1)^g\lambda_g/\zeta(1-2g)$. ( Here $\zeta(s)$ denotes the Riemann 
zeta function.) 
\par
\smallskip
\noindent
{\sl Proof.} The intersection of the closure of ${\cal A}_{g-1,1}$
with the boundary divisor in $\tilde {\cal A}_g$ is the closure of
$B_g$. Using the rational equivalence of $B_g$ and ${\cal
A}_{g-1} \times \{j\}$ for generic $j$ and the fact that the class of ${\cal
A}_{g-1} \times \{j\}$ on ${\cal A}_g^*$ was a multiple of $q_*\lambda_g$
the result follows by pull back. To find the multiple we look at
characteristic zero and integrate the $Sp(2g, \R)$-invariant forms
representing the $\lambda_i$. $\square$
\bigskip
\noindent
{\sl Examples.} We have in cohomology:
\smallskip
$g=1$ \qquad $[B_1]= 12 \lambda_1$;
\smallskip
$g=2$ \qquad $[B_2]=  120 \lambda_2$.
\noindent
\bigskip
Consider the moduli space ${\cal A}_g^{\prime}= \tilde{\cal A}_g^{(1)}$
of  rank
$\leq 1$ degenerations. Let $D^0$ be the closed subset corresponding to
rank 1 degenerations. The divisor $D^0$ has a morphism to 
$\phi:D^0 \to {\cal A}_{g-1}$ which  exhibits  $D^0$ as the universal
abelian variety over ${\cal A}_{g-1}$.  The fibre over $x \in {\cal
A}_{g-1}$ is the dual ${\hat X}_{g-1}$ of the abelian variety $X$ 
corresponding to $x$. The `universal' semi-abelian variety $G$ over
${\hat {\cal X}} _{g-1}$ is the $\G_m$-bundle obtained from the Poincar\'e
bundle $P \to {\cal X}_{g-1} \times {\hat {\cal X} }_{g-1}$ by deleting
the zero-section. We have the maps
$$
G=P -\{ (0)\} \to {\cal X}_{g-1} \times \hat{\cal X}_{g-1} 
\to {\cal A}_{g-1}.
$$
\smallskip
The divisor $D^0$ contains the subvariety $B_{g}$ corresponding to the
trivial extensions 
$$
1\to \G_m \to G \to X_{g-1}\to 0
$$
and this is a codimension $g$ cycle $B_{g} \sim {\cal A}_{g-1}$ in ${\cal
A}_g^{\prime}$.
\smallskip
Consider now the cotangent bundle to $G$ at the zero section $t$ of 
$G \to {\hat X}_{g-1}$. We have an exact sequence
$$
0 \to q^*\E_{g-1} \to {\E_g}_{|D^0} \to U\to 0,
$$
with $U$ a pull back of a line bundle on ${\cal A}_{g-1}^*$. Now $U$ is trivial
since the restriction of $\E_g$ to $B_{g}$ is a direct sum of $\E_{g-1}$ and a
trivial line bundle. 
\smallskip
Consider  ${\cal X}'$,  the compactified family of semi-abelian varieties
over the moduli stack ${\cal A}_g^{\prime}$ of rank $\leq 1$
degenerations. A semi-abelian variety which is a $\G_m$-bundle over an
abelian variety $X_{g-1}$ is compactified by taking the $\P^1$-bundle
associated to the $\G_m$-bundle and then identifying the $0$- and the
$\infty$-section under a shift $\xi \in X_{g-1}$. The image of the zero
section of the $\P^1$-bundle maps to a codimension $2$ cycle $\Delta$,
the locus of singular points of the fibres of ${\cal X}_g^{\prime} \to
{\cal A}_g^{\prime}$. 
We have 
$$
\Delta \cong {\cal X}_{g-1} \times_{{\cal A}_{g-1}}{\hat{\cal X}}_{g-1}.
$$
We analyze $\Omega^1= \Omega_{{\cal X}'_g / {\cal A}_g'}$. We have an exact
sequence $$
0 \to \Omega^1 \to \pi^*(\E)  \to  {\cal F} \to 0,\eqno(5)
$$
Here ${\cal F}$ is a sheaf with support on $\Delta$, the codimension $2$ cycle.
Let $u$ be a fibre coordinate on a $\G_m$ bundle over ${\cal X}_{g-1}$. A
section of
$\pi^*(\E)$ is given locally by $du/u$. Pull a section back to
the $\P^1$-bundle and take the residue along the $0$- and $\infty$-section. 
The residue map yields an isomorphism on ${\cal A}_g^{\prime}$
$$
{\cal F} \cong O_{\tilde{\Delta}},
$$
where $\tilde{\Delta}$ is the etale double cover of $\Delta$ corresponding to
choosing the branches ($0$~and~$\infty$). 
\bigskip
\noindent
\proclaim (1.11) Main Theorem. The tautological subring of $\tilde{\cal
A}_g$ in $CH_{\Q}^*(\tilde{\cal A}_g)$ is isomorphic to~$R_g$. \par
\bigskip
In order to prove this  we shall apply  GRR to the $\Theta$-divisor
again. We can do this  on a level cover of $\tilde{\cal
A}_g$ for a line bundle $L=O(T)$ trivialized along the zero
section. But we start in codimension $1$ and therefore we work on ${\cal
A}_g^{\prime}$. There we have:
$$
ch(\pi_!(L^{\otimes n}))= \pi_*(e^{nT}\cdot
{\rm Td}^{\vee}(O_{\Delta})^{-1}){\rm Td}^{\vee}(\E). \eqno(6)
$$
In particular, for $n=1$ we have
$$
ch(\pi_!(L))= \pi_*(e^{T}\cdot
{\rm Td}^{\vee}(O_{\Delta})^{-1}){\rm Td}^{\vee}(\E).\eqno(7)
$$
We write
$$
ch(\pi_!(L^{\otimes n})) = 1 + \theta_1^{(n)} + \theta_2^{(n)}+ \ldots  
$$
and set $\theta_1^{(1)}=\theta$. In equation (7)  we compare terms of
codimension $1$:
$$
 \theta_1^{(n)}= -{\lambda_1\over 2} \cdot \pi_*[e^{nT} \cdot {\rm
Td}^{\vee}(O_{\Delta})^{-1})]_0 + \pi_*[e^{nT}\cdot {\rm
Td}^{\vee}(O_{\Delta})^{-1})]_1 
$$
\smallskip
\noindent
\proclaim (1.12)  Lemma. We have $\theta=-\lambda_1/2+ \delta/8$, 
where $\delta$ is the class of the `boundary' $D$. \par
\smallskip
\noindent
{\sl Proof. } First we do the case $g=1$. We have
$$
1+ \theta = (\pi_*(T + T^2/2) + \pi_*(\Delta/12))(1-\lambda_1/2).
$$
Let $S$ be the zero-section. By Kodaira's results on elliptic surfaces
we have $S^2=-\lambda_1$; so the normalized $T$ is $T= S +
\pi^*(\lambda_1)$ with
$T^2= S\cdot
\pi^*(\lambda_1)$, i.e. $\pi_*(T^2/2)=\lambda_1/2$.  We find
$$
\theta =-\lambda_1/2+ \lambda_1/2 + \delta/12.
$$  
We can rewrite this as
$$
\theta =  -\lambda_1/2 + \delta/8.
$$
\smallskip
For general $g$ we have a priori  $\theta= -\lambda_1/2 +a\delta $. Restriction 
to the space of products of elliptic curves ${\cal A}_{1,\ldots,1}$ gives
$a=1/8$.
$\square$
\smallskip
\noindent
As a corollary we find 
$$
\theta_1^{(n)}= n^g(-\lambda_1/2) + (n^{g+1} + 2 n^{g-1})\delta/24.
$$
\noindent
To finish the proof we now work on the whole space $\tilde{\cal A}_g$
and we also introduce a (suitable) level
$\ell$ structure and then apply Grothendieck-Riemann-Roch on the moduli
space
$\tilde{\cal A}_g[\ell]$. The Hodge bundle is a pull back, but since the
natural maps
$\tilde{\cal A}_g[\ell] \to \tilde{\cal A}_g$ are ramified over the
divisor $D$ one can separate the contributions from the `interior' and the
boundary in the analogue of (6) by their dependence on $\ell$. This 
leads to the fact that the pull back of the formula $e^{-\lambda_1/2}=
{\rm Td}^{\vee}(\E)$ holds on the spaces $\tilde{\cal A}_g[\ell]$. We
refer to [EFG] for the details of the proof.
\bigskip
A top monomial $u^{\alpha}$ in $R_g$ can be written as a multiple
$m_{\alpha} u_1u_2\ldots u_g$ with $m_{\alpha}\in \Z $ using the relation
(1). Therefore, the degree $\deg \lambda^{\alpha}$
equals $m_{\alpha } \deg \lambda_1\lambda_2 \ldots
\lambda_g([\tilde{\cal A}_g])$. The ring $R_g$ is also the Chow ring
(and cohomology ring) of the Lagrangian Grassmann variety $Y_g$ of
maximal isotropic subspaces of a symplectic complex vector space of
dimension $2g$, see [P] and the references there. 
If we identify $R_g$ with $CH_{\Q}^*(Y_g)$ then we
find as a corollary:
\bigskip
\noindent
\proclaim (1.13) Theorem.  (Proportionality Principle of Hirzebruch-Mumford)
The characteristic numbers of the Hodge bundle are proportional to those of the
tautological bundle on $Y_g$:
$$
\lambda^{\alpha}([\tilde{\cal A}_g])= (-1)^G {1 \over 2^g}
\prod_{k=1}^g \zeta(1-2k) \cdot u^{\alpha}([Y_g]) 
$$
with $G= g(g+1)/2$.\par
\smallskip
\noindent
{\sl Proof.} We have $\lambda^{\alpha}([\tilde{\cal A}_g])=
c(g)\times u^{\alpha}([Y_g]$ for some constant $c(g)$. To evaluate the
constant one can compare the Euler number of ${\cal A}_g$ (in a
suitable sense; actually for $\Omega(\log D)= {\rm Sym}^2(\E)$)  and the
Euler number of $Y_g$. The Euler number of $Y_g$ equals $2^g$. We know the
Euler number of ${\cal A}_g$ by the  work of Siegel and Harder:
$$
\eqalign{\chi(Sp(2g, \Z))&= {\# W_{Sp(\R)} \over \# W_{U(\R)}} 
\prod_{j=1}^g {\zeta(1-2j)\over 2}
= {2^g g!\over g!} \prod_{j=1}^g {\zeta(1-2j)\over 2}\cr
&=\zeta(-1)\zeta(-3)\cdots \zeta(-2g+1).\cr}
$$
This proves the result. $\square$
\smallskip
Define a proportionality factor by 
$$
p(g)= (-1)^G \prod_{j=1}^g {\zeta(1-2j)\over 2}.
$$
Since $\deg u_1u_2\ldots u_g=1$
we find $\deg \lambda_1 \lambda_2 \cdots \lambda_g = p(g) $.
Similarly, using the structure of $R_g$ we have
$$
{\lambda_1^G \over \prod_{k=1}^g
\zeta(1-2k)}=(-1)^{G} {u_1^G\over 2^g}.
$$ 
and thus get 
$$
\lambda_1^G=p(g)  G!\prod_{k=1}^g{1 \over (2k-1)!!}. 
$$
This has to be interpreted with care (in the orbifold
sense). 
\vskip 1.0pc \noindent 
{\sl Some examples.} We have
$$
\eqalign{p(0) = 1 &\cr
p(1) = 1/24 \qquad &\deg \lambda_1=1/24,\cr
p(2) = 1/5760 \qquad &\deg \lambda_1^3=1/2880,\cr
p(3)  = 1/2903040 \qquad &\deg \lambda_1^6=1/181440.\cr}
$$
For the classical Proportionality principle of Hirzebruch we refer to his
collected works, [Hi 1]. In a letter to Atiyah Hirzebruch sketches
how one obtains a copy of $R_g$ in the cohomology of suitable
compact quotients of the Siegel upper half space, cf.\ [Hi~2].
\bigskip 
Proposition (1.2)
shows that the top Chern class $\lambda_g$ of $\E$ vanishes in the rational
Chow group $CH_{\Q}^g({\cal A}_g)$. So $\lambda_g$ is a torsion class on
${\cal A}_g$. Mumford proved in [M1] that the order of $\lambda_1$ in
$CH^1({\cal A}_1)$ is $12$. In [EFG] we shall prove the following bound
for the order of the torsion class $\lambda_g$. 
\smallskip
\proclaim (1.14) Definition.  Let $n_g$ be the greatest common divisor of
all $p^{2g}-1$ where $p$ runs through all primes greater than $2g+1$.
 \par
\smallskip
We have a little lemma.
\smallskip
\proclaim (1.15) Lemma. We have 
$$ 
\prod_{i=1}^g n_i = \prod_{p\hbox{  prime} }
([{2gp\over p-1}]!)_p \quad  
(= \hbox{\rm multiple of denominator of } p(g) \, ).
$$ \par
\smallskip
\proclaim (1.16) Theorem. On ${\cal A}_g$ we have: 
$(g-1)!\left( \prod_{i=1}^gn_i \right) \lambda_g=0$. \par \smallskip
\smallskip
\noindent
{\bf Example}
i) For $g=1$ we get $24\lambda_1=0$ which is off by a factor 2.
ii) For $g=2$ we get $24\cdot(16\cdot3\cdot5)\lambda_2=0$.
iii) For $g=3$ we get
$2\cdot24\cdot(16\cdot3\cdot5)\cdot(8\cdot9\cdot7)\lambda_3=0$.
\smallskip
We refer to [EFG] for other cycle relations in the Chow ring of
$\tilde{\cal A}_g$.
\bigskip
\eject
\noindent
{\bf \S 2. The Cycle Classes of the Ekedahl-Oort  Strata.} 
\bigskip
\noindent
Ekedahl and Oort introduced a stratification of the moduli space
${\cal A}_g \otimes \F_p$, cf. Oort's paper [O3] in this volume.  It
is defined by analyzing for abelian varieties $X$ the action of Frobenius
and Verschiebung on the group scheme
$X[p]$, the kernel of multiplication by $p$. The
strata include the well-known loci $V_f$ of abelian varieties of $p$-rank
$\leq f$  (for $0 \leq f \leq g$) and the loci 
$T_a$ of abelian varieties with $a$-number $\geq a$ (for 
$0\leq a \leq g$).
\smallskip
We shall describe these strata in a somewhat different way using the Hodge
bundle $\E$. We then can apply theorems of Porteous type to calculate the
cycle classes of these loci. Such Porteous type formulas are obtained by
applying results of Fulton and of Pragacz on degeneracy maps between
symplectic bundles. The cycle classes all lie in the tautological ring.
We get explicit formulas for the loci $V_f$ and $T_a$ generalizing 
Deuring's formula for the number of supersingular elliptic curves and
the formula for the number of superspecial abelian varieties  (with
$a=g$).
\smallskip
In the following we shall study the moduli space
${\cal A}_g \otimes \F_p$  of principally polarized abelian varieties in
characteristic $p$. For simplicity we shall write ${\cal A}_g$ instead of ${\cal
A}_g \otimes \F_p$.
\smallskip
Recall the {\sl canonical filtration} on the de Rham cohomology of a
principally polarized abelian variety $X$ as defined by Ekedahl and Oort.
We write
$G= H_{dR}^1(X)$ on which we have a $\sigma$-homomorphism $F$ and a
$\sigma^{-1}$-homomorphism $V$. For a moment we shall ignore
the $\sigma^{\pm}$-linearity. We have $FV=VF=0$.  The
$2g$-dimensional space is provided with a symplectic form
$\langle \, , \,\rangle$ and $F$ and $V$ are adjoints:
$\langle Vg, g' \rangle= \langle g, Fg' \rangle$. For any
subspace $H$ of $G$ we have $(VH)^{\bot} = F^{-1}(H^{\bot})$.
The spaces $V(G)= F^{-1}(0)$ and $F(G)=V^{-1}(0)$ are
maximally isotropic subspaces of dimension $g$. 
\smallskip
To construct the canonical filtration one starts with $0
\subset G$ and constructs finer filtrations by adding the
images of $V$ and the orthogonal complements of the spaces
present. This process stops. The filtration obtained is stable under $V$
and under ${}^{\bot}$, hence under $F^{-1}$ as well and is
called the {\sl canonical filtration}. The canonical filtration 
$$
0 \subset C_1 \subset ... \subset C_r \subset C_{r+1}
\subset ... \subset C_{2r} 
$$
satisfies $ C_r= V(C_{2r})$ and $C_{r-i}^{\bot} = C_{r+i}$. This
filtration can be refined to a so-called {\sl final filtration\/ } by 
choosing a $V$- and ${}^{\bot}$-stable filtration of length
$2g$ which refines the canonical one. In general, there is no  unique
choice for a final filtration. We thus get a  filtration
$$
0 \subset G_1 \subset ... \subset G_g \subset G_{g+1}
\subset ... \subset G_{2g} 
$$
which satisfies  $V(G_{2g})=G_g$,  $G_{g-i}^{\bot}=G_{g+i}$ and now also
$\dim(G_i)=i$. The associated {\sl final type} is the
increasing and surjective map
$$
\nu: \{ 0,1,2, \ldots ,2g\} \to \{ 0,1,2, \ldots, g\}
$$
satisfying
$$
\nu(2g-i)= \nu(i)-i+g \qquad {\rm for} \quad 0 \leq i \leq g\eqno(9)
$$
obtained by $\nu(i)=\dim(V(G_i))$ and $\nu(0)=0$. The {\sl canonical
type} is the restriction of
$\nu$ to the integers which arise as dimensions of the $C_i$. Although
the final filtration is not unique, the final type is. (The dimensions
between two steps of the canonical filtration either remain constant
or grow each step by $1$.) 
\smallskip
\noindent
{\bf Example}. Let $X$ be an abelian variety with $p$-rank $f$ and
$a(X)=1$ (equivalently, on $ G_g$ the operator $V$ has rank $g-1$ and
semi-simple rank $g-f$). Then the canonical type is given by the numbers
$\{ {\rm rank} (C_i)\}$
$$
 = \{ 0<f<f+1<\ldots <g-1<g<g+1 < \ldots
<2g-f-1<2g-f<2g\}
$$
and
$$
\eqalign{ \nu(f)=f, \, & \, \nu(f+1)=f, \cr 
\nu(f+2)=f+1, \ldots, \nu(g)=g-1, & \ldots, \nu(2g-f-1)=g-1,
\cr
\nu(2g-f)=g, \,  &  \, \nu(2g)=g.\cr}
$$
\bigskip
\noindent
It is not
difficult to see that there is a bijection between the set of final
types and the set of canonical types and we have
$2^g$ of them. Note that in view of (9) the function $\nu$ is determined
by its restriction to $\{1,2,\ldots, g\}$. This restriction is again
denoted by $\nu$. 
\smallskip
The Ekedahl-Oort stratification of ${\cal A}_g$ is obtained
by looking for each geometric point of ${\cal A}_g$ what the
canonical type (or final type) is. The set $Z_{\nu}$ of all
abelian varieties which have given final type $\nu$ is
locally closed and these $Z_{\nu}$ define a stratification, cf. [O3]. 
\smallskip
I prefer to describe the combinatorial datum of a final type  $\nu$ by a 
{\sl  partition\/} $\mu=\{g\geq \mu_1 > \mu_2 > \ldots > \mu_r \}$ as
follows:
$$
\mu_j = \# \{ i : 1 \leq i \leq g , \nu(i)\leq i-j \};
$$
equivalently, we can visualize it by the associated Young-type diagram
with
$\mu_j$ squares in the
$j$-th layer (i.e. by putting a stack of $i-\nu(i)$ squares in
position $i$):
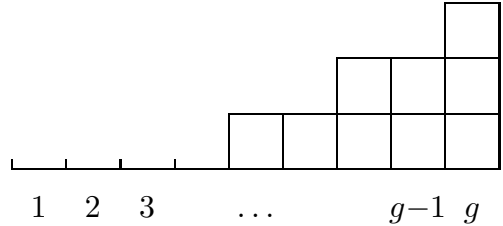
$$
\def\mt{6mm}
\def\mtstreepje{1mm}
\setbox0=\vbox{%
	\hrule width\mt
	\hbox to\mt{\vrule height\mt\hfil}}%
\setbox1=\hbox to\wd0{\hfil}%
\setbox2=\hbox to\wd0{\vrule height\mtstreepje\hfil}%
\def\tbox#1{\hbox to\wd0{\hss$#1$\hss}}%
\vbox{%
\hbox{\copy1\copy1\copy1\copy1\copy1\copy1\copy1\copy1\copy0\vrule}
\nointerlineskip
\hbox{\copy1\copy1\copy1\copy1\copy1\copy1\copy0\copy0\copy0\vrule}
\nointerlineskip
\hbox{\copy2\copy2\copy2\copy2\copy0\copy0\copy0\copy0\copy0\vrule}
\nointerlineskip
\hrule\vskip3mm
\hbox{\tbox1\tbox2\tbox3\copy1\tbox{\ldots}\copy1\copy1\tbox{g{-}1}\tbox{g}}}%
$$
This example  corresponds to 
$\{ \nu(i) : i=1,...,g\} = \{1,2,\ldots,g-5, g-5, g-4, g-4, g-3, g-3 \}$ 
and to  $\mu= \{5 , 3, 1\}$.
\smallskip
\noindent
{\bf (2.1) Definition.} We call a partition $\mu$ {\sl admissible} if
$g\geq \mu_1 > \mu_2 > \ldots >\mu_r >0$. 
We call the function $\nu \colon \{1,\ldots , g\} \to \{
1,\ldots, g\}$ {\sl admissible} if
$$
0\leq \nu (i) \leq \nu (i+1) \leq
\nu(i)+1\leq g+1
$$ for $1 \leq i \leq g$.
The number $|\mu|:= \sum_i \mu_i$ is called the {\sl area
of the diagram}. 
\bigskip
\noindent
The  notions `admissible' diagram  (or partition)
and `admissible function' $\nu$ and final type are all equivalent.
The set of admissible diagrams carries a partial ordering in the obvious
way: $\mu \geq \mu'$ if $\mu_i \geq \mu_i'$ for all $i$.
\smallskip 
We shall now give another approach to these strata
by defining them globally on a flag space over ${\cal A}_g$. Our starting
point is the observation that  if we define for a final
filtration $E_i$ ($1\leq i \leq g$) another filtration by $F_g= \ker (V)=
V^{-1}(0)$, by  $F_{g+i}= V^{-1}(E_i)$ for $i=1,...,g$ and
by $F_{g-i}= F_{g+i}^{\bot}$ then we have
$$
V(E_i)\subseteq E_{\nu(i)} \Longleftrightarrow 
E_i \subseteq V^{-1}(E_{\nu(i)}) = F_{g+\nu(i)} 
\Longleftrightarrow \dim (E_i \cap F_{g+\nu(i)}) \geq i.\eqno(10)
$$
\smallskip
Working now globally, we  let $S$ be a scheme in characteristic
$p$ and  let ${\cal X}\to S$ be an abelian variety over $S$ with principal
polarization. Then we consider the de Rham cohomology sheaf 
${\cal H}_{dR}^1({\cal X}/S)$. It is defined as the hyper-direct image $ 
{\cal R}^1\pi_*(O_{\cal X} \to \Omega_{{\cal X}/S}^1)$. 
It is a locally free sheaf of rank $2g$ on $S$  equipped with a
non-degenerate alternating form (cf.\ [O])
$$
\langle \, , \, \rangle \colon {\cal H}_{dR}^1({\cal X}/S) \times {\cal
H}_{dR}^1({\cal X}/S) \to O_S.
$$
Indeed, the polarization (locally in the \'etale topology given by a
relatively ample line bundle on ${\cal X}/S$) provides us with a
symmetric homomorphism $\rho \colon {\cal X} \to \hat {\cal X}$ and the
Poincar\'e bundle defines a perfect pairing between $ {\cal
H}_{dR}^1({\cal X}/S)
$ and ${\cal H}_{dR}^1(\hat{\cal X}/S)$. Moreover, we have an exact
sequence of locally free sheaves $$
0 \to \pi_*(\Omega^1) \to {\cal H}_{dR}^1({\cal X}/S) 
\to R^1\pi_*O_{\cal X} \to 0.
$$
We shall write $\H$ for the sheaf ${\cal H}_{dR}^1({\cal X}_g /
{\cal A}_g)$ and ${\cal X}$ for ${\cal X}_g$. We thus have an
exact sequence
$$
0 \to \E \to \H \to \E^{\vee} \to 0.
$$
The relative Frobenius  $F: {\cal X} \to {\cal X}^{(p)}$ and
the Verschiebung $V:{\cal X}^{(p)} \to {\cal X}$ satisfy $F\cdot V =
p\cdot {\rm id}_{{\cal X}^{(p)}}$ and $V \cdot F= p \cdot {\rm id}_{\cal
X}$ and they induce maps in cohomology, again denoted by $F$ and $V$: 
$$
F:\H^{(p)} \to \H \qquad {\hbox {\rm and } } \quad V: \H \to \H^{(p)}.
$$
Of course, we have  $FV=0$ and $VF=0$ and  $F$ and $V$ are adjoints.
This implies that $ {\rm Im}(F)= \ker (V)$ and ${\rm Im}(V) = \ker (F) $
are maximally isotropic subbundles of $\H$ and $\H^{(p)}$. Moreover, 
since $dF=0$ on ${\rm Lie}({\cal X})$ it follows  that
$F=0$ on $\E$ and thus ${\rm Im}(V) = \E^{(p)}$.
Verschiebung thus provides us with a bundle map (again denoted by $V$):
$V: \H \to \E^{(p)}$.
\smallskip
Consider the space of symplectic flags ${\cal F}= {\rm Flag}(\H)$ on
the  bundle $\H$. This space is fibred by the spaces  ${\cal F}^{(i)}$
of partial flags
$$
\E_i \subsetneq \E_{i+1} \subsetneq \ldots \subsetneq \E_g.
$$
So ${\cal F}^{(1)}= {\rm Flag}(\H)$ and ${\cal F}^{(g)}= {\cal A}_g$ and
there are natural maps
$$
\pi_{i,i+1}\colon {\cal F}^{(i)} \to {\cal F}^{(i+1)}.
$$
The fibres are Grassmann varieties of dimension $i$. The space ${\cal
F}^i$ is equipped with a universal flag. On ${\cal F}$ the Chern classes 
of the bundle $\E$ decompose into their roots: 
$$
\lambda_i= \sigma_i(\ell_1,\ldots, \ell_g)\qquad {\rm with} \qquad
\ell_i = c_1(\E_{i}/ \E_{i-1}).
$$
\smallskip
Given  an arbitrary flag of subbundles
$$
0 \subsetneq \E_1 \subsetneq\ldots \subsetneq \E_g=\E \eqno(11)
$$
with rank$(\E_i)=i$ we can extend this to a symplectic filtration on
$\H$  by putting
$$
\E_{g+i}= (\E_{g-i})^{\bot}.
$$
By base change we can transport this filtration to $\H^{(p)}$.
\smallskip
We introduce a second filtration
by starting with the isotropic subbundle
$$
\F_g=\ker (V)= V^{-1}(0) \subset \H
$$
and continuing with
$$
\F_{g+i}= V^{-1}(\E_i^{(p)}) \qquad {\rm for} \quad 1\leq i \leq g.
$$
We extend it to a symplectic filtration by setting $\F_{g-i}=
(\F_{g+i})^{\bot}$. We thus have two filtrations $\E_{\bullet}$ and
${\F}_{\bullet}$ on $\H$. 
\smallskip
\noindent
{\bf Example.} i) Let $X$ be an ordinary abelian variety. Then ${\rm
Lie}(X)= {\rm Lie}(\mu)$ with $\mu$ the multiplicative subgroup scheme 
of $X[p]$ of rank $p^g$. It follows that  $V$ is  invertible on ${\rm
Lie}(X)^{\vee}= \omega(X)$, i.e. $\F_g \cap \E_g= (0)$.
ii) If $X$ is a superspecial abelian variety (i.e. $X$ without its
polarization is a product of supersingular elliptic curves) then $V=0$ on
$\omega(X)$ so that ${\rm rk} (\E_i \cap \F_g)= i$. 
\smallskip
These two (extreme) examples show that the respective position of the
two filtrations $\E_{\bullet}$ and $\F_{\bullet}$  for an abelian variety
$X$ gives information on the structure of the kernel of multiplication by
$p$ on $X$. These respective positions are encoded by  a combinatorial
datum, e.g.\ $\nu$ or to be more precise, by  an element of a Weyl group.
We shall associate strata to such data.
\smallskip
To either $\nu$ or $\mu$  we now associate an element of the Weyl
group of the symplectic group. The Weyl group $W_g$ of type $C_g$ in
Cartan's terminology  is isomorphic to the semi-direct product $S_g
\ltimes   (\Z / 2\Z)^g$, where $S_g$ acts on $(\Z / 2\Z)^g$ by
permuting the $g$ factors. Another description of this group  is as the 
subgroup of $S_{2g}$ of elements which map any symmetric
$2$-element subset of the form $\{ i , 2g+1-i \}$ of 
$\{ 1, \ldots, 2g\}$ to a subset of the same type:
$$
W_g = \{ \sigma \in S_{2g} \colon \sigma(i) + \sigma(2g+1-i) = 2g+1 
{ \hbox { \rm for } } i=1,\ldots, g \}.
$$
An element in this Weyl group has a {\sl length} and a {\sl
codimension}:
$$
\ell(w) = \#\{i < j : w(i) > w(j) \} + \# \{ i < j :
w(i)+w(j)> 2g+1\}
$$
and
$$
{\rm codim}(w) =\#\{i < j : w(i) < w(j) \} + \# \{ i < j :
w(i)+w(j) < 2g+1\}.
$$
We have the equality
$$
\ell(w)+ {\rm codim}(w) = g^2.
$$
\bigskip
To a function $\nu $ we associate the following  element of the Weyl
group, a the permutation of $\{ 1,2, \ldots, 2g\}$ : let 
$$
S= \{ i_1,i_2,\ldots \}= \{ 1 \leq i \leq g : \nu(i)=\nu(i-1)\}
$$ 
with $i_1 < i_2 < \ldots $ given  in increasing order. Let 
$$
S^c= \{ j_1, j_2, \ldots \}
$$
be the elements of $\{ 1,2, \ldots, g\}$ not in $S$, in increasing
order. Then one gets the permutation of $S_{2g}$ defined by $\nu$ by
writing $g+k$ at position $j_k$ for $k=1,2,\ldots$ and $k$ at position
$i_k$ for $k=1,2,\ldots$. We finish off by putting $2g+1-k$ at
position $2g+1-i$ if $k$ is written at position $i$. 
We obtain a sequence $s$ which is a permutation of $\{1,2,\ldots, 2g\}$.
\smallskip 
Alternatively, using diagrams, we can describe the element $w=w_{\mu}$ as
follows. Let $t_i$ be the operator on the set of  diagrams which is
`remove the top box of the $i$-th column'. We let the $t_i$ act from the
right. Given an admissible diagram $\mu$ we consider the complementary
diagram  $$
\mu^c= \{ g,g-1,\ldots, 1\} - \mu
$$
which is also an admissible diagram (partition).  We can  successively
apply operators $t_i$ to it  such that after every step we obtain an
admissible diagram and such that after  $|\mu^c|$ steps we obtain
the empty  diagram $(\mu^c) \cdot t_{i_1 }\cdots t_{i_{\ell}} =
\emptyset$. We remove first the top layer, then the next one and so
on. For $1\leq i < g$ let $s_i\in S_{2g}$ be the permutation
$(i, i+1)(2g-i, 2g+1-i)$ and let $s_g=(g, g+1) \in S_{2g}$. Now
associate to the diagram
$\mu$ the element of the Weyl group:
$$
\mu \mapsto w_{\mu}= s_{i_1} \cdots s_{i_{\ell}}.
$$
Each diagram thus yields an element in the Weyl
group. The admissible partitions yield  $2^g$ elements of the $2^gg!$
elements of $W_g$.
\bigskip
\noindent
{\bf Example}. $g=3$
\noindent
\settabs\+&\indent&$\mu $ xxxxxxxxx & $\nu $ xxxxxxxx & $
\{1,2,3,4,5,6\}$ maps to xx & xx$\ell$ xxx & element xxxxxxxx\cr
\+&&$ { } \mu $ & $ { } \nu $  & $ [1,2,3,4,5,6] $ maps to & $\ell$ &
$w_{\mu}$
\cr
\bigskip
\+&&$\emptyset$ & $\{ 1, 2, 3\} $ & $[4,5,6,1,2,3]$ & $6$ &
$s_3s_2s_3s_1s_2s_3$\cr
\+&&$\{ 1 \} $ & $\{ 1, 2, 2 \}$ & $[4,5,1,6,2,3]$ & $5$ &
$s_2s_3s_1s_2s_3$\cr 
\smallskip
\+ && $\{ 2 \} $ & $\{ 1, 1, 2 \}$ & $[4,1,5,2,6,3]$ & $4$ &
$s_3s_1s_2s_3$\cr
\smallskip
\+ && $\{ 3 \} $ & $\{ 0, 1, 2 \}$ & $[1,4,5,2,3,6]$ & $3$ &
$s_3s_2s_3$\cr
\smallskip
\+ && $\{ 2,1 \} $ & $\{ 1,1,1 \}$ & $[4,1,2,5,6,3 ]$ & $3$ &
$s_1s_2s_3$\cr
\smallskip
\+ && $\{ 3 \} $ & $\{ 0, 1, 2 \}$ & $[1,4,5,2,3,6]$ & $3$ &
$s_3s_2s_3$\cr
\smallskip
\+ && $\{ 3, 1  \} $ & $\{ 0, 1, 1 \}$ & $[1,4,2,5,3,6]$ & $2$ &
$s_2s_3$\cr
\smallskip
\+ && $\{ 3,2  \} $ & $\{ 0, 0,1 \}$ & $[1,2,4,3,5,6]$ & $1$ &
$s_3$\cr
\smallskip
\+ && $\{ 3, 2, 1  \} $ & $\{ 0, 0, 0  \}$ & $[1,2,3,4,5,6]$ & $0$ &
$1$\cr
\bigskip
\bigskip
We can associate to the map $V: \H \to \E^{(p)}$ and an element $w$ 
a degeneracy locus ${\cal U}_{w}$ and in particular to an
admissible diagram $\mu$ a degeneracy locus ${\cal U}_{\mu}={\cal
U}_{w_{\mu}}$  in ${\cal F}={\rm Flag}(\H)$ . Intuitively, ${\cal
U}_w$ is defined as the locus of points $x$ such that at $x$ we have 
$$
\dim (\E_i \cap \F_j) \geq \# \{ a \leq i: w(a) \leq j \} \quad {\hbox
{\rm for all } } \quad 1 \leq i,j \leq g
$$
or equivalently, that ${\rm ker}(V)\cap \E_i \geq i-\nu(i)$ (= the number
of squares in the diagram $\mu$ in position $i$). For the precise
definition we refer to Fulton [F]. Note that for the admissible 
diagrams we do not use the full filtration of $\F$, but only $\F_g$.
\bigskip
We first look at the case of the empty diagram $\mu= \emptyset$ or
equivalently, that $\nu=\{1,2,3,\ldots, g\}$.  The degeneracy conditions
say that (cf. (10))
$$
V(\E_i) \subseteq \E_i^{(p)} \qquad i=1,\ldots, g.
$$
i.e. we are looking at the space ${\cal U}_{\{ \emptyset\} }$ of symplectic
filtrations on $\E$ which are compatible with the action of $V$. The
codimension of this space in ${\cal F}$ is  $g(g-1)/2$, hence $\dim
(U_{\{ \emptyset \} }) = g(g+1)/2= \dim ({\cal A}_g)$. We have a
finite map $\pi : U_{\{ \emptyset \}}\to {\cal A}_g$ of degree
$$
\deg (\pi) = (1+p)(1+p+p^2)\ldots (1+p+p^2+\ldots + p^{g-1}).
$$
This space ${\cal U}= U_{ \{\emptyset \} }$ can be seen as a component of
the moduli space $\Gamma_0(p)$.  Indeed, a filtration of the subgroup
scheme $X[p]$ corresponds to a filtration on $H_{dR}^1$. Fix $g$ and
let $\Gamma_0(p)^{\prime}$ be the functor that associates to a
scheme $S$ the set of isomorphism classes of principally polarized abelian
schemes $X$ over $S$  plus a $V$-stable symplectic filtration on $X[p]$. 
\smallskip
\noindent
\proclaim (2.2) Proposition. The degeneracy cycle ${\cal U}_{\emptyset}$
is the algebraic stack representing the functor  $\Gamma_0(p)^{\prime}$.
The stack ${\cal U}_{\emptyset}$ is fibred by finite morphisms
$$
{\cal U}_{\emptyset}= {\cal U}^{(1)} {\buildrel\pi_{1,2}
\over\lllongrightarrow} \,{ \cal U}^{(2)} {\buildrel\pi_{2,3} \over
\lllongrightarrow} \ldots {\buildrel \pi_{g-1,g} \over
\lllongrightarrow } \, { \cal U}^{(g)}= {\cal A}_g.
$$
with $\deg(\pi_{i,i+1})= 1+p+\ldots + p^{i}$. It contains the
degeracy loci ${\cal U}_{\mu}$ for all $2^g$ partitions $\mu$. \par
\smallskip 
The stacks ${\cal U}^{(i)}$ come
with universal (partial) flags $\E_i \subsetneq \ldots \subsetneq \E_g$.
We denote by $\lambda_j(i)= c_j(\E_i)$ the Chern class in
$CH_{\Q}^j({\cal U}^{(i)})$. We also have tautological quotient bundles
$L_i=\E_i/ \E_{i-1}$. We denote by $l_i$ the Chern class of $L_i$. We
have  $(\lambda_j)(i)= \sigma_j(l_1,\ldots, l_i)$, the $j$-th
elementary symmetric function of the $l_1, \ldots, l_i$. 
\bigskip
Next we look at the cases where $\mu$ is a partition of the form
$\mu =\{ \mu_1= g-f \}$. The corresponding degeneracy loci classify
the loci of $p$-rank $\leq f$. It is well-known by Oort (cf. [O-N]) that
the codimension of $V_f$  in ${\cal A}_g$ equals $g-f$.  They admit the
following explicit description.  The pullback of
$V_{g-1}$ to
${\cal U}_{\emptyset}$ consists of $g$ components, say $Z_1, ...,Z_g$. An
abelian variety of
$p$-rank $g-1$ and $a=1$ has a unique subgroup scheme $\alpha_p$. The
index $i$ of $Z_i$ indicates where for the generic point of $Z_i$ the
$\alpha_p$ can be found: $\alpha_p \subset G_i$, $\alpha_p \not\subset
G_{i-1}$.
\smallskip
Then the pull back of  $V_{f}$ consists of  $$ \pi^{-1}(V_f)=
\sum_{\#S = g-f} Z^S, $$
where for a subset $S \subset \{1, \ldots, g \}$ the cycle $Z^S$ is
defined as
$$
Z^S= \cap_{i\in S} Z_i.
$$
Then one sees easily that the cycle $V_{g-f}$ on ${\cal A}_g$  is
obtained as the push forward of $Z_g\cap \ldots \cap Z_{g-f+1}$.
\proclaim (2.3) Lemma.
The cycle class of $Z_i$ is equal to $(p-1)\ell_i$. \par
\noindent
{\sl Proof.} It is described as the locus where $\det (V) : L_i \to
L_i^{(p)}$ vanishes. By viewing $\det(V)$ as a section of $L_i^{(p)}
\otimes L_i^{-1}$ the result follows. $\square$
\smallskip
We find after a calculation:
$$
\pi^*([V_f]) = (p-1)^{g-f} \lambda_{g-f} 
$$ 
on ${\cal U}_{\emptyset}$. Using the push forward of the $\pi_{i,
i+1}$ on the classes $\ell_i$ and extending the result to a
compactification
$\tilde{\cal A}_g$ we get: 
\smallskip
\noindent
\proclaim (2.4) Theorem.  The cycle class of $V_f$, the $p$-rank $\leq
f$ locus, in the Chow ring $CH_{\Q}^*(\tilde{\cal A}_g)$,  is given by
$[V_f] = (p-1)(p^2-1)\ldots (p^{g-f}-1) \lambda_{g-f}$. \par
\noindent
\proclaim (2.5) Corollary. (Deuring Mass Formula) Let $g=1$. We have
$$
\sum_E {1 \over \# {\rm Aut}(E) } = {p-1\over 24},
$$
where the sum is over the isomorphism classes over ${\bar{\F}}_p$ of
supersingular elliptic curves.
\par
\noindent
{\sl Proof.} The formula gives $(p-1)\lambda_1$ and by (1.10) this
equals $\delta/12$. The class of $\delta$ is equivalent to $1/2$
times the class of a `physical' point of the $j$-line because the
degenerate elliptic curve corresponding to $\delta$ has $2$
automorphisms. $\square$
\bigskip
Another case where we can find an explicit formula is the case of the
locus $T_a$. This corresponds to the case $\mu= \{ a, a-1, \ldots,
2,1\}$. But here we can work directly on ${\cal A}_g$. 
The locus $T_a$ on
${\cal A}_g$ may be defined as the locus $$ \{ x \in {\cal A}_g : {\rm
rank}(V)|_{\E_g} \leq g-a \}. $$
We have $T_{a+1} \subset T_a$ and $\dim (T_g)=0$. 
\smallskip
Pragacz and Ratajski, cf.\  [P-R],  have developed
formulas for the degeneracy locus for the rank of a self-adjoint bundle
map of symplectic bundles globalizing the results in isotropic Schubert 
calculus from [P]. Before we apply  their result to our case we
have to introduce some notation.
\smallskip
Define for a vectorbundle $A$ with Chern classes
$a_i$ the expression
$$
Q_{ij}(A) := a_ia_j + 2\sum_{k=1}^j (-1)^k a_{i+k}a_{j-k}
\quad { \rm for } \quad i>j.
$$
For a admissible  partition $\beta  = (\beta_1, ... , \beta_r)$
(with $r$ even, $\beta_r$ may be zero) we set
$$
Q_{\beta} = {\rm Pfaffian}(x_{ij}),
$$ 
where the $(x_{ij})$ is an anti-symmetric matrix with
$ x_{ij} = Q_{\beta_i,\beta_j} $. Applying the Pragacz-Ratajski formula
to  our  situation gives the following result:
\smallskip
\proclaim (2.6) Theorem. The class of the reduced locus $T_a$ of abelian 
varieties with $a$-number $\geq a$  is given by 
$$
\sum Q_{\beta}(\E^{(p)})\cdot Q_{\rho(a) -\beta}(\E^*) ,
$$
where the sum is over the admissible partitions $\beta$ contained in
the partition $\rho(a) = (a, a-1, a-2, ..., 1)$. \par
\noindent
{\bf Example:} 
$$ 
\eqalign{ 	[T_1]=& p\lambda_1 - \lambda_1\cr
	[T_2] =& (p-1)(p^2+1) (\lambda_1 \lambda_2) - (p^3-1) 2 \lambda_3 \cr
& ... \cr
	[T_g]  =& (p-1)(p^2+1) \ldots (p^g+ (-1)^g) \lambda_1 \lambda_2
\cdots \lambda_g.\cr }
$$
As a corollary we find a classical result of Ekedahl (cf.\ [E]) on the 
number of principally polarized abelian varieties with $a=g$:
\noindent
\proclaim (2.7) Corollary. We have
$$
\sum_X {1 \over \#{\rm Aut}(X)} = (-1)^G 2^{-g} \big[\prod_{j=1}^g 
(p^j+(-1)^j)\big] \cdot \zeta(-1) \zeta(-3) \ldots \zeta(1-2g),
$$
where the sum is over the isomorphism classes (over ${\bar{\F}}_p$) of
principally polarized abelian varieties of dimension $g$ with $a=g$. \par
\smallskip
\noindent
{\sl Proof.} Combine the formula for $T_g$ with the
Proportionality Theorem. $\square$ 
\bigskip
For each element $w$ we now find a degeneracy locus ${\cal U}_w$. In
particular we find such a locus  ${\cal U}_{\mu}$ for each partition
$\mu$ in ${\cal F}$ and the ${\cal U}_{\mu}$  actually lie in ${\cal
U}_{\emptyset}$.  It is known that $T_g$ is zero-dimensional, cf. (2.7).
This implies that the codimension of each ${\cal U}_{\mu}$ equals ${\rm
codim}(w(\mu))= |\mu|$. We can apply a theorem of Fulton to determine the
class of the degeneracy locus $U_{\mu}$ in $CH_{\Q}^*({\cal F})$.
For an admissible diagram  $\mu= \{ \mu_1> \ldots > \mu_r>0\}$ 
Fulton  defines a determinant (`Schur function') 
$$
\Delta_{\mu}(x_i)= \det ( x_{\mu_i+j-i})_{1\leq i,j\leq r};
$$
this is a polynomial with integer coefficients in the commuting
variables $x_1, x_2, \ldots$. We define a `double Schubert function' by 
putting
$$
\Delta(x,y) = \Delta_{g,g-1,\ldots, 1}\big( \sigma_i(x_1,\ldots,
x_g)+\sigma_i(y_1,\ldots, y_g)\big).
$$
The operators $\partial_i$ (`divided difference operators') on the
ring $\Z[x_1,x_2, \ldots, ]$ are defined by setting for  $F(x) \in
\Z[x_1,\ldots,x_g]$ 
$$
\cases{ \partial_i(F(x)) = { F(x)- F(s_i(x))\over x_i-x_{i+1}} & if
$i<g$\cr
\partial_g(F(x))= {F(x)-F(s_g(x))\over 2x_g} & if $i=g$.\cr}
$$
We put
$$
\Delta= \Delta(x,y) = \Delta_{g,g-1,\ldots, 1}\big( \sigma_i(x_1,\ldots,
x_g)+\sigma_i(y_1,\ldots, y_g)\big).
$$
We shall apply the result of Fulton and state an abstract formula for our
cycle classes. In principle one then can calculate the push forward
algorithmically. But it seems difficult to get closed formulas for the
cycle classes of these push forwards. 
\smallskip
By applying a theorem of Fulton [F] we find:
\smallskip
\noindent
\proclaim (2.8) Theorem. Let $\mu$ be an admissible diagram whose
corresponding element $w_{\mu }$ in the Weyl group is written as
$s_{i_{\ell}}\cdots s_{i_1}$. The cycle class of the degeneracy locus
$U_{\mu }$ in $CH_{\Q}^*({\cal F})$ is given by
$$
u_{\mu }=\partial_{i_\ell}\cdots \partial_{i_1}\big( (\prod_{i+j\leq
g}(x_i-y_{i}))\cdot \Delta(x,y) \big)_{\big| \{ x_i=pl_i, y_i=-l_i \} }
$$
\par
\proclaim (2.9) Corollary. The push forward of the class of $[{\cal
U}_{\mu}]$ under $\pi$ is given by $\pi_*(u_{\mu})$ and is  a multiple
of the class of the reduced cycle $Z_{\mu}$, the multiple being equal
to the number of final filtrations refining  the canonical
filtration associated to $\mu$. The class belongs
to the tautological ring. 
\par 
\bigskip 
We can calculate these classes
from this formula in an algorithmic way. But is seems difficult to get
closed formulas in the general case. We list here the formulas for $g=3$.
The multiplicity is given by the factor in square brackets.
\smallskip
\noindent
{\bf Formulas for $g=3$}.
$$
\pi_{*}(\{ \emptyset \} )=[(1+p)(1+p+p^2)]
$$
$$
   \pi_{*}(\{ 1\} )= [\left( 1 + p \right)]\times  \, \left( 1 - p
\right)
     \lambda_{1}$$
$$\pi_{*}(\{ 2\} )=
    (p-1)(p^2-1)\,\lambda_{2}$$
$$
   \pi_{*}(\{ 1,2\} )=[(1+p)] \times \{ \, ( 1 - p ) \,
( 1 + {{p}^2} ) \,\lambda_{1}\,\lambda_{2} -
     2\, ( -1 + p^3 ) \,
      \lambda_{3} \}$$
$$
   \pi_{*}(\{ 3\} )= (p-1)(p^2-1)(p^3-1)  \,\lambda_{3}$$
$$
   \pi_{*}(\{ 1,3\} )=[( 1 + p ) ] \times {{( -1 + p
) }^2}\, \,
     ( 1 - p + {{p}^2} ) \,\lambda_{1}\,\lambda_{3}$$
$$
   \pi_{*}(\{ 2,3\} )={{(-1 + p) }^3}\, (1 + p) \,
     ( -1 + p - {{p}^2}) \, ( 1 + p + {{p}^2}) \,
     \lambda_{2}\,\lambda_{3}$$
$$
   \pi_{*}(\{ 1,2,3\} )=[(1+p^3)] \times (p-1) \, (p^2+1)
    (p^3-1)  \,\lambda_{1}\,\lambda_{2}\,\lambda_{3}.
$$
\bigskip
The canonical filtration on $H_{dR}^1(X)$  is the most
economical one with respect to the operation $V$. The types of this
filtration correspond to the elements of $(\Z / 2\Z)^g$ in the Weyl group
$(\Z / 2\Z)^g \rtimes S_g$.  The other types of filtrations  are
obtained by applying an element of $S_g$ to the right. Therefore, in the
closure of a stratum ${\cal U}_{\mu}$ we find also ${\cal U}_w$ with $w
\not\in (\Z / 2\Z)^g$. This explains the phenomenon observed by Oort in
[O3], 4.6. We can use Pieri-type formulas to get results on the boundary
of strata. In particular one obtains the result of
Ekedahl-Oort that the class $\lambda_1$ is torsion on the open strata.
We refer to [EFG].
\smallskip
Some strata can be described in detail. Consider the locus
associated to the partition
$$
\mu = \{ g, g-1,  \dots, 3,2\}.
$$
This is the penultimate stratum and it is of dimension $1$. It is
proved by Ekedahl and Oort that this locus is connected. 
\smallskip
Define $h(g)$ by $h(1)=1$ and $h(g)=h(g-1)\times{p^g+(-1)^g\over
p+(-1)^g}$. Then we have
$$
\pi_*(U_{[g,g-1,...,1]})=h(g)[T_g]
$$
with $[T_g]$ given as above. 
\bigskip
\noindent
\proclaim (2.10) Theorem. The cycle class of $\pi_*({\cal U}_{\mu})$ for
$\mu=\{g,g-1,\ldots, 3,2\}$ is given by  
$$
{(p-1)\over (p^2+1)} \times h(g) \times \left(\prod_{i=1}^g (p^i+(-1)^i)
\right) \times \lambda_2 \lambda_3 \ldots \lambda_g.
$$
In ${\cal U}_{\emptyset}$ this locus consists of a configuration
of $\P^1$'s.
\par
\smallskip
Since the degree of $\lambda_1$ on each copy of $\P^1$ is $p-1$ one can
compute the number of components. This locus is highly reducible. But
some loci are irreducible:
\bigskip
\noindent
\proclaim (2.11) Theorem. For $a<g$ the locus $T_a$ is irreducible. In
particular, the locus $T_1= V_{g-1}$ is irreducible. \par
\noindent
{\sl Proof.} By Ekedahl-Oort (cf. [O3]) we know that the locus
corresponding to the diagram $\{g,g-1, \ldots, 3,2\}$ is connected. This
implies that $T_a$ is connected for $a<g$. We know that ${\rm Sing}(T_a)
\subseteq T_{a+1}$, hence of codimension $>1$. Actually, one can
describe the normal bundle to $T_a-T_{a+1}$.  By a theorem of Hartshorne
(cf. [H]) this implies that $T_a- {\rm Sing}(T_a)$ is connected.
$\square$
\smallskip
\noindent
The irreducibility of $V_{g-1}$ was also observed by Oort, cf. [O3].
\bigskip
\noindent
{\bf \S 3 Some additional results.}
\bigskip
\noindent
We describe as an example the canonical type for hyperelliptic curves
of $2$-rank $0$ in characteristic $p=2$.
\smallskip
\noindent 
\proclaim (3.1) Lemma. A hyperelliptic curve  $C$ of genus $g$  
and of $2$-rank $0$ over
$ k = {\bar k}$ of characteristic $p=2$  can be written as
$$
y^2+y=xP(x^2)=x(a_1x^2+a_2x^4+...+a_{g-1}x^{2g-2}+x^{2g}).\eqno(13)
$$ 
\par 
\smallskip
\noindent
\proclaim (3.2) Theorem. For a hyperelliptic curve as in (13)
the canonical type is described by the corresponding partition 
$\mu= [g,g-2,g-4,...]$. In particular, the canonical filtration is
independent of the coefficients $a_i$. 
\par
\smallskip
We refer for the proof to [EFG].
\bigskip
For $g=1$ and $g=2$ the supersingular locus
coincides with $V_0$ and thus occurs as a stratum in the Ekedahl-Oort
stratification. For $g\geq 3$ this is no longer true. The definition of
supersingular locus involves  also the higher filtrations, i.e. it uses  
not only $X[p]$, but the higher group schemes $X[p^i]$ as well. But in
a number of cases we can determine the class of this locus. We state
here the result  for the case $g=3$. \smallskip 
\noindent 
\proclaim (3.3) Theorem. The class of the supersingular
locus $S_3$  in ${{\cal A}}_3^*\otimes \F_p$ is  
$$
[S_3]= (p-1)(p^2-1)(p^3-1)(p-1)(p^2+1) \lambda_1 \lambda_3.
$$
\par
The proof uses the explicit description by Li and Oort of the moduli of
principally polarized supersingular abelian varieties, cf. [L-O]. We
refer the reader to [EFG] for the proof and formulas for other genera.
\bigskip
\noindent
{\bf References}
\smallskip
\noindent
[B-S] A.\ Borel, J-P. Serre: Le th\'eor\`eme de Riemann-Roch (d'apr\`es
Grothendieck). {\sl Bull.\ Soc.\ Math.\ France \bf 86} (1958), 97-136.
\smallskip
\noindent
[E] T.\ Ekedahl: On supersingular curves and abelian varieties.
{\sl Math.\ Scand.\ \bf 60} (1987), 151-178.
\smallskip
\noindent
[EFG] T.\ Ekedahl, C.\ Faber, G.\ van der Geer: Manuscripts in
preparation
\smallskip \noindent
[E-O] T.\ Ekedahl, F.\ Oort :
Connected subsets of a moduli space of abelian varieties. Preliminary
version of a preprint.
\smallskip
\noindent
[F-C] G.\ Faltings, C-L.\ Chai: Degeneration of abelian varieties.
Ergebnisse der Math. 22. Springer Verlag 1990.
\smallskip
\noindent
[F] W.\ Fulton: Determinantal formulas for orthogonal and
symplectic degeneracy loci. To appear in {\sl J.\ Diff.\ Geometry} 
(1996).
\smallskip
\noindent
[F-L]
W.\ Fulton, R.\ Lazarsfeld: On the connectedness of degeneracy loci and
special divisors. {\sl Acta Mathematica \bf 146}, (1981), p.\ 271-283.
\smallskip
\noindent
[H] R.\ Hartshorne: Complete intersections and connectedness. {\sl
Amer.\ J.\ Math.\ \bf 84} (1962), 497-508.
\smallskip
\noindent
[Hi 1] F.\ Hirzebruch: Automorphe Formen und der Satz von Riemann-Roch.
In {\sl Symposium Internacional de Topologia Algebrica} (M\'exico
1956), 129-144. M\'exico: La Universidad Nacional Aut\'onoma de
M\'exico 1958. (= Collected Papers I , Springer Verlag,  p. 345)
\smallskip
\noindent
[Hi 2] F.\ Hirzebruch: Kommentar zu `Elliptische Differentialoperatoren
auf Mannigfaltigkeiten'. In: Collected Papers II, Springer Verlag,  p.
773.
\smallskip
\noindent
[L-O] K-Z.\  Li,  F.\ Oort:
Moduli of supersingular abelian varieties.Preprint University Utrecht,
Nr. 824 (1993).
\smallskip
\noindent
[M-B] L.\ Moret-Bailly: Pinceaux de vari\'et\'es ab\'eliennes. 
Ast\'erisque 129 (1985).
\smallskip
\noindent
[M1] D.\ Mumford: Picard groups of moduli problems. In {\sl
Arithmetical Algebraic Geometry}. Ed.\ O.F.G.\ Schilling.
Harper and Row, 1965, p.\ 33-81. 
\smallskip
\noindent
[M2] D.\ Mumford: Hirzebruch's Proportionality Theorem in the
non-compact case. {\sl Inv.\ Math.\ \bf 42} (1977) , 239-272.
\smallskip
\noindent
[M3] D.\ Mumford: On the Kodaira dimension of the Siegel modular
variety. In: Algebraic Geometry-Open Problems. SLNM 997, 348-375.
\smallskip
\noindent
[N-O] P.\ Norman, F. Oort:  Moduli of abelian varieties. {\sl Ann.\
Math.\  \bf 112}  (1980), 413-439.
\smallskip
\noindent
[O] T.\ Oda : The first de Rham cohomology group and Dieudonn\'e
modules. Ann.\ Sci.\ ENS (1969), p. 63-135.
\smallskip
\noindent
[O1] F.\ Oort: Subvarieties of moduli spaces. {\sl Inv.\ Math.\ \bf 24}
(1974), 95-119.
\smallskip
\noindent
[O2] F.\ Oort: Complete subvarieties of moduli spaces. In: {\sl Abelian
Varieties} (W.\ Barth, K.\ Hulek, H.\ Lange, eds.), de Gruyter Verlag,
Berlin, 1995, p. 225-235.
\smallskip
\noindent
[O3] F.\ Oort: A stratification of a moduli space of polarized
abelian varieties in positive characteristic. Preprint May 1996.
\smallskip
\noindent
[P] P.\ Pragacz: Algebro-geometric applications of Schur S- and Q-
polynomials. In: Topics in Invariant Theory. S\'eminaire d'Alg\`ebre
Dubreil-Malliavin 1989-1990. SLNM  1478, 130-191 (1991).
\smallskip
\noindent
[P-R] P.\ Pragacz, J.\ Ratajski :Formulas for Lagrangian and 
orthogonal degeneracy loci; The $\tilde Q$-polynomials approach. To
appear in {\sl Comp.\ Math.\ }
\bigskip
Gerard van der Geer
\smallskip\indent
Faculteit WINS
\smallskip\indent
Universiteit van Amsterdam
\smallskip\indent
Plantage Muidergracht 24
\smallskip\indent
1018 TV Amsterdam
\smallskip\indent
The Netherlands
\smallskip\indent
e-mail: {\tt geer@fwi.uva.nl}
\end